\documentclass[11pt]{article}
\usepackage{CJK}
\usepackage{amsmath}
\usepackage{amsmath,amssymb,latexsym,color}
\usepackage[mathscr]{eucal}
\newtheorem{thm}{Theorem}[section]
\newtheorem{ob}[thm]{Observation}
\newtheorem{defin}[thm]{Definition}
\newtheorem{prop}[thm]{Proposition}
\newtheorem{lemma}[thm]{Lemma}
\newtheorem{cor}[thm]{Corollary}

\newtheorem{example}[thm]{Example}

\newcommand{\proof}{{\it Proof.\quad}}
\newcommand{\qed}{\hfill\Box\medskip}
\newcommand{\rdim}{{\rm rdim}}
\begin{document}

\title{\bf On the metric dimension and  fractional metric dimension for hierarchical product of graphs}
\author{Min Feng\quad Kaishun Wang\footnote{Corresponding author. E-mail address: wangks@bnu.edu.cn}\\
{\footnotesize   \em  Sch. Math. Sci. {\rm \&} Lab. Math. Com. Sys.,
Beijing Normal University, Beijing, 100875,  China}}
\date{}
\maketitle

\begin{abstract}

A set of vertices $W$ {\em resolves} a graph $G$ if every vertex of
$G$ is uniquely determined by its vector of distances to the
vertices in $W$. The {\em metric dimension} for $G$, denoted by
$\dim(G)$, is the minimum cardinality of a resolving set of $G$. In
order to study the metric dimension  for the hierarchical product
$G_2^{u_2}\sqcap G_1^{u_1}$ of two rooted graphs $G_2^{u_2}$ and $
G_1^{u_1}$, we first introduce a new parameter, the {\em rooted
metric dimension} $\rdim(G_1^{u_1})$  for a rooted graph
$G_1^{u_1}$. If $G_1$ is not a path with an end-vertex $u_1$, we
show that $\dim(G_2^{u_2}\sqcap
G_1^{u_1})=|V(G_2)|\cdot\rdim(G_1^{u_1})$, where $|V(G_2)|$ is the
order of $G_2$. If $G_1$ is a path with an end-vertex $u_1$, we
obtain some tight inequalities for $\dim(G_2^{u_2}\sqcap
G_1^{u_1})$.  Finally, we show that similar results hold for the
fractional metric dimension.



\end{abstract}
\bigskip

\bigskip

\section{Introduction}
All graphs considered in this paper are nontrivial and connected.
For a graph $G$, we often denote by $V(G)$ and $E(G)$ the vertex set
and the edge set of $G$, respectively. For any two vertices $u$ and
$v$ of $G$, denote by $d_G(u,v)$ the distance between $u$ and $v$ in
$G$, and write $R_G\{u,v\}=\{w\mid w\in V(G),d_G(u,w)\neq
d_G(v,w)\}$. If the graph $G$ is clear from the context, the
notations $d_G(u,v)$ and $R_G\{u,v\}$ will be written $d(u,v)$ and
$R\{u,v\}$, respectively. A subset $W$ of $V(G)$ is  a {\em
resolving set} of $G$ if $W\cap R\{u,v\}\neq\emptyset$ for any two
distinct vertices $u$ and $v$. A {\em metric basis} of $G$ is a
resolving set of $G$ with minimum cardinality. The cardinality of a
metric basis of $G$ is the {\em metric dimension}  for $G$, denoted
by $\dim(G)$.

Metric dimension was introduced independently by Harary and Melter
\cite{Ha}, and by Slater \cite{Sl1}. As a graph parameter it has
numerous applications, among them are computer science and robotics
\cite{Kh}, network discovery and verification \cite{Bee}, strategies
for the Mastermind game \cite{Chv} and combinatorial optimization
\cite{Se}. Metric dimension has been heavily studied, see \cite{Ba}
for a number of references on this topic.

 The problem of finding
the metric dimension  for a graph was formulated as an integer
programming problem independently by Chartrand et al. \cite{cha},
and by Currie and Oellermann \cite{cu}. In graph theory,
fractionalization of integer-valued graph theoretic concepts is an
interesting area of research (see \cite{sc}). Currie and Oellermann
\cite{cu} and   Fehr et al. \cite{fe} defined fractional metric
dimension as the optimal solution of the linear relaxation of the
integer programming problem. Arumugam and Mathew \cite{Ar} initiated
the study of the fractional metric dimension  for graphs. For more
information, see \cite{Aru,fl,fw}.

Let $g: V(G)\longrightarrow[0,1]$ be a real value function. For
$W\subseteq V(G)$, denote $g(W)=\sum_{v\in W}g(v)$. The {\em weight}
of $g$ is defined by $|g|=g(V(G))$. We call $g$ a {\em resolving
function} of $G$ if $g(R\{u,v\})\geq 1$ for any two distinct
vertices $u$ and $v$. The minimum weight of a  resolving function of
$G$  is called  the {\em fractional metric dimension}  for $G$,
denoted by $\dim_f(G)$.

 It was noted in \cite[p.
204]{Ga} and \cite{Kh} that determining the metric dimension  for a
graph is an NP-complete problem. So it is  desirable to reduce the
computation for the  metric dimension  for product graphs  to the
computation for some parameters of the factor graphs; see \cite{Ca2}
for cartesian products, \cite{ja} for lexicographic products,
and \cite{ye} for corona products. Recently, the fractional
metric dimension  for the above three products was studied in
\cite{Aru,fl,fw}.

In order to model some real-life complex networks, Barri\`ere et al.
\cite{Ba1} introduced the hierarchical product of graphs and showed
that it is associative. A {\em rooted graph} $G^u$ is the graph $G$
in which one vertex $u$, called {\em root vertex}, is labeled in a
special way to distinguish it from other vertices. Let $G_1^{u_1}$
and $G_2^{u_2}$ be two rooted graphs. The {\em hierarchical product}
$G_2^{u_2}\sqcap G_1^{u_1}$ is the rooted graph with the vertex set
$\{x_2x_1\mid x_i\in V(G_i), i=1,2\}$, having the root vertex
$u_2u_1$, where $x_2x_1$  is adjacent to $y_2y_1$ whenever $x_2=y_2$
and $\{x_1,y_1\}\in E(G_1)$, or $x_1=y_1=u_1$ and $\{x_2,y_2\}\in
E(G_2)$.  See \cite{ju,no,rab,ras} for more information.

 In this
paper, we study the (fractional) metric dimension for the
hierarchical product $G_2^{u_2}\sqcap G_1^{u_1}$ of rooted graphs
$G_2^{u_2}$ and $G_1^{u_1}$. In Section 2, we introduce a new
parameter, the rooted metric dimension $\rdim(G^{u})$ for a rooted
graph $G^{u}$. If $G_1$ is not a path with an end-vertex $u_1$, we
show that $\dim(G_2^{u_2}\sqcap
G_1^{u_1})=|V(G_2)|\cdot\rdim(G_1^{u_1})$. If $G_1$ is a path with
an end-vertex $u_1$, we obtain some tight inequalities for
$\dim(G_2^{u_2}\sqcap G_1^{u_1})$. In Section 3, we show that
similar results hold for the fractional metric dimension.

\section{Metric dimension}

In order to study the  metric dimension for the hierarchical product
of graphs, we first introduce  the rooted metric dimension  for a
rooted graph.

 A {\em rooted resolving set} of a rooted graph $G^u$
is a subset $W$ of $V(G)$ such that $W\cup\{u\}$ is a resolving set
of $G$. A {\em rooted metric basis} of $G^u$ is a  rooted resolving
set of $G^u$ with the minimum cardinality. The cardinality of a
rooted metric basis of $G^u$ is called {\em rooted metric dimension}
for  $G^u$, denoted by $\rdim(G^u)$. The following observation
is obvious.

\begin{ob}\label{ob1}
If there exists a metric basis of $G$  containing $u$, then
$\rdim(G^u)=\dim(G)-1.$ If any metric basis of $G$ does not contain
$u$, then $\rdim(G^u)=\dim(G)$.
\end{ob}

For graphs $H_1$ and $H_2$ we use $H_1\cup H_2$ to denote the
disjoint union of $H_1$ and $H_2$ and $H_1+H_2$ to denote the graph obtained from the disjoint
union of $H_1$ and $H_2$ by joining every vertex of $H_1$ with every vertex of $H_2$.

\begin{ob}\label{ob2}
Let $G$ be a graph of order $n$. Then $1\leq\dim(G)\leq n-1.$
Moreover,

{\rm(i)} $\dim(G)=1$ if and only if $G$ is the path $P_n$ of length $n$.

{\rm(ii)} $\dim(G)=n-1$ if and only if $G$ is the complete graph $K_n$ on $n$ vertices.
\end{ob}

\begin{prop}\label{cha}
{\rm\cite[Theorem 4]{cha}} Let $G$ be a graph of order $n\geq 4$.
Then $\dim(G)=n-2$ if and only if $G=K_{s,t}$ $(s,t\geq 1)$,
$G=K_s+\overline{K_t}$ $(s\geq 1,t\geq 2)$, or $G=K_s+(K_1\cup K_t)$
$(s,t\geq 1)$, where $\overline{K_t}$ is a null graph and $K_{s,t}$
is a complete bipartite graph.
\end{prop}
\begin{prop}\label{path}
Let $G^u$ be a rooted graph of order $n$. Then $0\leq\rdim(G^u)\leq
n-2$. Moreover,

{\rm (i)} $\rdim(G^u)=0$ if and only if $G=P_n$ and $u$ is one of its
end-vertices.

{\rm (ii)} $\rdim(G^u)=n-2$ if and only if $G=K_n$,
or $G=K_{1,n-1}$and $u$ is the centre.
\end{prop}
\proof If $G$ is a complete graph, by Observation~\ref{ob2} (ii) we have
$\dim(G)=n-1$ . Observation~\ref{ob1} implies that $\rdim(G^u)=n-2$.
If $G$ is not a complete graph, then $1\leq\dim(G)\leq n-2$, which
implies that $0\leq\rdim(G^u)\leq n-2$ by Observation~\ref{ob1}.

(i) Since $\rdim(G^u)=0$ if and only if $\{u\}$ is a metric basis of
$G$, by Observation~\ref{ob2} (i), (i) holds.

(ii) Suppose $\rdim(G^u)=n-2$. Then $\dim(G)=n-1$ or $n-2$. If
$\dim(G)=n-1$, then $G=K_n$. Now we consider
$\dim(G)=n-2$. If $n=3$, then $\dim(G)=1$, which implies that $G=K_{1,2}$
and $u$ is the centre. Now suppose $n\geq
4$. Then $G$ is one of graphs listed in Proposition \ref{cha}. If
$s,t\geq 2$ or $G=K_s+(K_1\cup K_t)$, then there exists a metric
basis containing $u$, which implies that $\rdim(G^u)=n-3$, a
contradiction. Hence $G=K_{1,n-1}$. Since any
metric basis of $K_{1,n-1}$ does not contain the centre, the vertex
$u$ is the centre of $K_{1,n-1}$. The converse is routine. $\qed$

Next, we study the  metric dimension for the hierarchical product of
graphs.

 Let $G_1^{u_1}$ and $G_2^{u_2}$ be two rooted graphs. For
any two vertices $x_{2}x_{1}$ and $y_{2}y_{1}$ of $G_2^{u_2}\sqcap
G_1^{u_1}$,   observe that
\begin{equation}\label{1}
d(x_{2}x_{1},y_{2}y_{1})=\left\{
\begin{array}{ll}
d_{G_1}(x_1,y_1),                 &\textup{if}~x_{2}=y_{2},\\
d_{G_2}(x_{2},y_{2})+d_{G_1}(x_1,u_1)+d_{G_1}(y_1,u_1),&\textup{if}~x_{2}\neq y_{2}.
\end{array}\right.
\end{equation}

\begin{lemma}\label{resolve}
Let $x_2x_1$ and $y_2y_1$ be two distinct vertices of
$G_2^{u_2}\sqcap G_1^{u_1}$.

{\rm (i)} If $x_{2}=y_{2}$, then
$$
R\{x_2x_1,y_2y_1\}=\left\{
\begin{array}{ll}
\{x_2z\mid z\in R_{G_1}\{x_1,y_1\}\}, &\textup{ if  }u_1\not\in R_{G_1}\{x_1,y_1\},\\
V(G_2^{u_2}\sqcap G_1^{u_1})\setminus\{x_2z\mid z\not\in R_{G_1}\{x_1,y_1\}\},&\textup{ if  } u_1\in R_{G_1}\{x_1,y_1\}.
\end{array}\right.
$$

{\rm(ii)} If $x_{2}\neq y_{2}$, then $\{x_2z,y_2z\}\cap R\{x_2x_1,y_2y_1\}\neq\emptyset$
for any $z\in V(G_1)$.
\end{lemma}
\proof (i) If $u_1\not\in R_{G_1}\{x_1,y_1\}$, then
$d_{G_1}(x_1,u_1)=d_{G_1}(y_1,u_1)$. By (\ref{1}), the inequality
$d(x_2x_1,z_2z_1 )\neq d(y_2y_1,z_2z_1)$ holds if and only if
$z_2=x_2$ and $d_{G_1}(x_1,z_1)\neq d_{G_1}(y_1,z_1). $ It follows
that $R\{x_2x_1,y_2y_1\}=\{x_2z\mid z\in
R_{G_1}\{x_1,y_1\}\}.$ If $ u_1\in R_{G_1}\{x_1,y_1\}$, then
$d_{G_1}(x_1,u_1)\neq d_{G_1}(y_1,u_1)$. By (\ref{1}), the equality
$d(x_2x_1,z_2z_1 )= d(y_2y_1,z_2z_1)$ holds if and only if $z_2=x_2$
and $d_{G_1}(x_1,z_1)= d_{G_1}(y_1,z_1).$ It follows that
$R\{x_2x_1,y_2y_1\}=V(G_2^{u_2}\sqcap G_1^{u_1})\setminus\{x_2z\mid
z\not\in R_{G_1}\{x_1,y_1\}\}.$

(ii) Suppose $x_2z\not\in R\{x_2x_1,y_2y_1\}$. Then
$d(x_2x_1,x_2z)=d(y_2y_1,x_2z)$. By (\ref{1}),
$$
d_{G_1}(x_1,z)=d_{G_2}(y_2,x_2)+d_{G_1}(y_1,u_1)+d_{G_1}(z,u_1)\geq
d_{G_2}(x_2,y_2)+d_{G_1}(y_1,z),
$$
which implies that
$$
d_{G_2}(x_2,y_2)+d_{G_1}(x_1,u_1)+d_{G_1}(z,u_1)\geq
2d_{G_2}(x_2,y_2)+d_{G_1}(y_1,z)>d(y_2y_1,y_2z).
$$
Hence, $y_2z\in
R\{x_2x_1,y_2y_1\}$, as desired. $\qed$

\begin{lemma}\label{metdimlow}
Let $G_1^{u_1}$ and $G_2^{u_2}$ be two rooted graphs.
Then
$$
\rdim(G_2^{u_2}\sqcap G_1^{u_1})\geq |V(G_2)|\cdot\rdim(G_1^{u_1}).
$$
\end{lemma}
\proof  Let $\overline W$ be a rooted metric basis
 of $G_2^{u_2}\sqcap G_1^{u_1}$. For $v\in V(G_2)$, write $\overline W_v=\{z\mid vz\in\overline
 W\}$.
For any two distinct vertices $x,y$ of $G_1$, there exists a vertex
$wz$ in $\overline W\cup\{u_2u_1\}$ such that $d(vx,wz)\neq
d(vy,wz).$ If $w=v$, by (\ref{1}) we get $d_{G_1}(x,z)\neq
d_{G_1}(y,z)$, which implies that $z\in (\overline
W_v\cup\{u_1\})\cap R_{G_1}\{x,y\}. $ If $w\neq v$,  by (\ref{1}) we
have $d_{G_1}(x,u_1)\neq d_{G_1}(y,u_1)$, which implies that $
u_1\in R_{G_1}\{x,y\}. $ Therefore, we have $(\overline
W_v\cup\{u_1\})\cap R_{G_1}\{x,y\}\neq\emptyset,$ which implies that
$\overline W_v$ is a rooted resolving set of $G_1^{u_1}$. Hence
$$
\rdim(G_2^{u_2}\sqcap G_1^{u_1})=|\overline W|=\sum_{v\in
V(G_2)}|\overline W_v|\geq |V(G_2)|\cdot\rdim(G_1^{u_1}),
$$
as desired.
$\qed$

\begin{thm}\label{metmain}
Let $G_1^{u_1}$ and $G_2^{u_2}$ be two rooted graphs. If $G_1$ is
not a path with an end-vertex $u_1$, then
$$
\dim(G_2^{u_2}\sqcap G_1^{u_1})=|V(G_2)|\cdot\rdim(G_1^{u_1}).
$$
\end{thm}
\proof By Lemma \ref{metdimlow}, we only need to prove that
\begin{equation}\label{7}
\dim(G_2^{u_2}\sqcap G_1^{u_1})\leq|V(G_2)|\cdot\rdim(G_1^{u_1}).
\end{equation}
 Let $W$ be a rooted metric basis of
$G_1^{u_1}$. Then $W\neq\emptyset$. Write $\overline W=\{vw\mid v\in
V(G_2),w\in W\}$. Note that $|\overline
W|=|V(G_2)|\cdot\rdim(G_1^{u_1})$. In order to prove (\ref{7}), we
only need to show that $\overline{W}$ is a resolving set of
$G_2^{u_2}\sqcap G_1^{u_1}$. It suffices to show that, for any two
distinct vertices $x_2x_1$ and $y_2y_1$ of $G_2^{u_2}\sqcap
G_1^{u_1}$,
\begin{equation}\label{5}
\overline{W}\cap R\{x_2x_1,y_2y_1\}\neq\emptyset.
\end{equation}
If $x_2=y_2$ and $u_1\not\in R_{G_1}\{x_1,y_1\}$, then $W\cap
R_{G_1}\{x_1,y_1\}\neq\emptyset$, by Lemma \ref{resolve} (i) we
obtain (\ref{5}). If $x_2=y_2$ and $u_1\in R_{G_1}\{x_1,y_1\}$, by
(\ref{1}) we have $vw\in\overline{W}\cap R\{x_2x_1,y_2y_1\}$ for any
$v\neq x_2$ and any $w\in W$, which implies that (\ref{5}) holds. If
$x_2\neq y_2$, then (\ref{5}) holds by Lemma \ref{resolve} (ii). We
accomplish our proof. $\qed$

Combining Observation~\ref{ob1} and Theorem~\ref{metmain}, we have
the following result.

\begin{cor}
Let $G_1^{u_1}$ and $G_2^{u_2}$ be two rooted graphs.

{\rm(i)} If there exists a metric basis of $G_1$ containing $u_1$ and $G_1$ is not a path, then
$$
\dim(G_2^{u_2}\sqcap G_1^{u_1})=|V(G_2)|(\dim(G_1)-1).
$$

{\rm(ii)} If any metric basis of $G_1$ does not contain $u_1$, then
$$
\dim(G_2^{u_2}\sqcap G_1^{u_1})=|V(G_2)|\dim(G_1).
$$
\end{cor}

The {\em binomial tree} $T_n$ is the hierarchical product of $n$
copies of the complete graph on two vertices, which is a useful data
structure in the context of algorithm analysis and designs
\cite{co}. It was proved that the metric dimension for a tree  can
be expressed in terms of its parameters in \cite{cha,Ha,Sl1}.

\begin{cor}\label{29}
Let $n\geq 2$. Then $\dim(T_n)=2^{n-2}$.
\end{cor}
\proof Note that $\dim(T_2)=1$. Now suppose $n\geq 3$. Since
$T_n=(K_2^0\sqcap\cdots\sqcap K_2^0)\sqcap (K_2^0\sqcap K_2^0)$ and
$\rdim(K_2^0\sqcap K_2^0)=1$, the desired result follows by Theorem
\ref{metmain}. $\qed$

We always assume that $0$ is one end-vertex of $P_n$. In the
remaining of this section, we shall prove some tight inequalities
for $\dim(G^{u}\sqcap P_n^0)$.

\begin{prop}\label{sqcappath}
Let $G^u$ be a rooted graph with diameter $d$.
Then
\begin{eqnarray}
&&\dim(G^u\sqcap P_n^0)\leq\dim(G^u\sqcap P_{n+1}^0) \textup{ for } 1\leq n\leq d-1,\label{8}\\
&&\dim(G^u\sqcap P_n^0)=\dim(G^u\sqcap P_{n+1}^0)\textup{ for } n\geq d.\label{9}
\end{eqnarray}
\end{prop}
\proof If $G=K_2$, then $G^u\sqcap P_n^0$ is the path, which implies
that (\ref{9}) holds. Now we only consider $|V(G)|\geq 3$. Suppose
that $\overline W_n$ is a metric basis of $G^u\sqcap P_{n}^0$. Let
$P_n=(z_0=0,z_1,\ldots,z_{n-1})$. Define $\pi_n:V(G^u\sqcap
P_{n+1}^0)\longrightarrow V(G^u\sqcap P_n^0)$ by
$$
\pi_n(vz_i)=\left\{
\begin{array}{ll}
vz_{n-1},&\textup{ if } i=n,\\
vz_i, &\textup{ if } i\leq n-1.
\end{array}\right.
$$
Then $\pi_n(\overline W_{n+1})$ is a resolving
set of $G^u\sqcap P_n^0$, which implies that $\dim(G^u\sqcap P_n^0)\leq\dim(G^u\sqcap P_{n+1}^0)$
for any positive integer $n$. So (\ref{8}) holds.

In order to prove (\ref{9}), we only need to show that $\overline W_n$ is a resolving
set of $G^u\sqcap P_{n+1}^0$ for $n\geq d$.
Pick any two distinct vertices $v_1z_i$ and
$v_2z_j$ of $G^u\sqcap P_{n+1}^0$. It suffices to prove that
\begin{equation}\label{12}
\overline W_n\cap R_{G^u\sqcap P_{n+1}^0}\{v_1z_i,v_2z_j\}\neq\emptyset.
\end{equation}

We claim that  there exist two distinct vertices $w_1$ and $w_2$ of
$G$ such that $\overline W_n\cap\{w_sz_k\mid 0\leq k\leq
n-1\}\neq\emptyset$ for $s\in\{1,2\}$. Suppose for the contradiction
that there exists a vertex $w\in V(G)$ such that $\overline
W_n\subseteq\{wz_k\mid 0\leq k\leq n-1\}$. If the degree of $w$ in
$G$ is one, then there exists a path $(w,x,y)$ in $G$. For any
$wz_k\in\overline W_n$, we have $d(xz_1,wz_k)=k+2=d(yz_0,wz_k)$,
contrary to the fact that $\overline W_n$ is a metric basis of
$G^u\sqcap P_{n}^0$. If the degree of $w$ in $G$ is at least two,
pick two distinct neighbors $x$ and $y$ of $w$ in $G$. Then
$d(xz_0,wz_k)=k+1=d(yz_0,wz_k)$ for any $wz_k\in\overline W_n$, a
contradiction. Hence our claim is valid.

Now we prove (\ref{12}). Without loss of generality, we may assume
that $0\leq i\leq j\leq n$. If $j\leq n-1$, then $R_{G^u\sqcap
P_{n+1}^0}\{v_1z_i,v_2z_j\}\supseteq R_{G^u\sqcap
P_{n}^0}\{v_1z_i,v_2z_j\}$; and so (\ref{12}) holds. Now suppose
$j=n$.

{\em Case 1.} $v_1=v_2$. By the claim, the set $\{w_1z_k\mid 0\leq
k\leq n-1\}$ or $\{w_2z_k\mid 0\leq k\leq n-1\}$ is a subset of
$R_{G^u\sqcap P_{n+1}^0}\{v_1z_i,v_1z_n\}$. So (\ref{12}) holds.

{\em Case 2.} $v_1\neq v_2$.

{\em Case 2.1.} $i=0$. By the claim, we can choose $w_sz_k\in\overline W_n$ with $w_s\neq v_2$. Then
$$
d(v_1z_0,w_sz_k)=d_G(v_1,w_s)+k\leq d+k\leq n+k<d_G(v_2,w_s)+n+k=d(v_2z_n,w_sz_k),
$$
which implies that
$w_sz_k\in R_{G^u\sqcap P_{n+1}^0}\{v_1z_0,v_2z_n\}$. So
(\ref{12}) holds.

{\em Case 2.2.} $i\geq 1$. Note that
$$
R_{G^u\sqcap P_{n+1}^0}\{v_1z_i,v_2z_n\}=R_{G^u\sqcap P_{n+1}^0}\{v_1z_{i-1},v_2z_{n-1}\}
\supseteq R_{G^u\sqcap P_{n}^0}\{v_1z_{i-1},v_2z_{n-1}\}.
$$
Then
(\ref{12}) holds. $\qed$

\begin{prop}\label{metpathbound}
For any rooted graph $G^u$, we have
\begin{equation}\label{6}
\dim(G)\leq\dim(G^u\sqcap P_n^0)\leq|V(G)|-1.
\end{equation}
\end{prop}
\proof Let $z$ be the other end-vertex of $P_n$. Fix a vertex
$v_0\in V(G)$ and write $\overline S=\{vz\mid v\in
V(G)\setminus\{v_0\}\}$. Since $\{z\}$ is a resolving set of $P_n$,
the set $\overline S$ resolves $G\sqcap P_n$ by (\ref{1}). Hence
$\dim(G^u\sqcap P_n^0)\leq |\overline S|=|V(G)|-1.$ Since $G^u$ is
isomorphic to $G^u\sqcap P_1^0$, Proposition \ref{sqcappath} implies
that $\dim(G)\leq\dim(G^u\sqcap P_n^0)$. $\qed$

For $m\geq 2$, we have $\dim(K_m^u\sqcap P_n^0)=m-1$. This shows
that the inequalities (\ref{8}) and (\ref{6}) are tight.

\begin{example}
For $m,n\geq 2$, we have $\dim(P_m^u\sqcap P_n^0)=2$.
\end{example}
\proof Write $P_k=(z_0=0,z_1,\ldots, z_{k-1})$. Then
$\{z_0z_{n-1},z_{m-1}z_{n-1}\}$ is a resolving set of $P_m^u\sqcap
P_n^0$. $\qed$

\begin{example}
Let $C_m$ be the cycle with length $m$. Then $\dim(C_m^u\sqcap
P_n^0)=2$.
\end{example}
\proof Let $P_n=(z_0=0,z_1,\ldots, z_{n-1})$ and
$C_m=(c_0,c_1,\ldots, c_{m-1},c_0)$. Then
$\{c_0z_{n-1},c_1z_{n-1}\}$ is a resolving set of $C_m^u\sqcap
P_n^0$. $\qed$

\section{Fractional metric dimension}

In order to study the fractional metric dimension for the
hierarchical product of graphs, we first introduce  the fractional
rooted metric dimension for a rooted graph.

Similar to the fractionalization of metric dimension,
we give a fractional version of the rooted metric dimension
 for a rooted graph.
Let $G^u$ be a rooted graph of order $n$. Write
$$
\mathcal{P}^u=\{\{v,w\}\mid v,w\in V(G), v\neq w, d(v,u)=d(w,u)\}.
$$
Suppose $\mathcal{P}^u\neq\emptyset.$ Write
$V(G)\setminus\{u\}=\{v_1,\ldots,v_{n-1}\}$ and $\mathcal
P^u=\{\alpha_1,\ldots,\alpha_m\}$. Let $A^u$ be the $m\times(n-1)$
matrix with
$$
(A^u)_{ij}=\left\{
\begin{array}{ll}
1,&\textup{ if }v_j \textup{ resolves } \alpha_i,\\
0,&\textup{ otherwise. }
\end{array}\right.
$$
The integer programming formulation of the rooted metric dimension
for $G^u$ is given by
\begin{eqnarray*}
&&\textup{Minimize }f(x_1,\ldots,x_{n-1})=x_1+\cdots+x_{n-1}\\
&&\textup{Subject to }A^u\mathbf{x}\geq\mathbf{1}
\end{eqnarray*}
where ${\mathbf x}=(x_1,\ldots,x_{n-1})^{\rm T},x_i\in\{0,1\}$ and
${\mathbf 1}$ is the $m\times 1$ column vector all of whose entries
are $1$. The optimal solution of the linear programming relaxation
of the above integer programming problem, where we replace
$x_i\in\{0,1\}$ by $x_i\in[0,1]$, gives the {\em fractional rooted
metric dimension} for $G^u$, which we denote by $\rdim_f(G^u)$.

Let $G^u$ be a rooted graph which is not a path with an end-vertex
$u$. A {\em rooted resolving function} of a rooted graph $G^u$ is a
real value function $g:V(G)\longrightarrow[0,1]$ such that
$g(R\{v,w\})\geq 1$ for each $\{v,w\}\in \mathcal{P}^u$. The {\em
fractional rooted metric dimension} for $G^u$ is the minimum weight
of a rooted resolving function of $G^u$.

\begin{prop}\label{root frac}
Let $G^u$ be a rooted graph which is not a path with an end-vertex
$u$. Then

{\rm(i)} $\rdim_f(G^u)\leq\rdim(G^u)$.

{\rm(ii)} $\rdim_f(G^u)\leq\frac{|V(G)|-1}{2}$.

{\rm(iii)} $\dim_f(G)-1\leq\rdim_f(G^u)\leq\dim_f(G)$.
\end{prop}
\proof (i)  Let $W$ be a rooted metric basis of $G^u$. Define
$g:V(G)\longrightarrow[0,1]$ by
$$
g(v)=\left\{
\begin{array}{ll}
1,&\textup{ if }v\in W,\\
0,&\textup{ if }v\not\in W.
\end{array}\right.
$$
For any $\{x,y\}\in\mathcal P^u$, there exists a vertex $v\in W$
such that $d(x,v)\neq d(y,v)$. Then $g(R\{x,y\})\geq g(v)=1$, which
implies that $g$ is a rooted resolving function of $G^u$. Hence
$\rdim_f(G^u)\leq |g|=|W|=\rdim(G^u)$.

(ii) The function
$g:V(G)\longrightarrow[0,1]$ defined by
$$
g(v)=\left\{
\begin{array}{ll}
0,&\textup{ if }v=u,\\
\frac{1}{2},&\textup{ if }v\neq u
\end{array}\right.
$$
is a rooted resolving function of $G^u$. Hence
$\rdim_f(G^u)\leq\frac{|V(G)|-1}{2}$.

(iii) It is clear that $\rdim_f(G^u)\leq\dim_f(G)$. Let $g$ be a
rooted resolving function of $G^u$. Then the function
$h:V(G)\longrightarrow[0,1]$ defined by
$$
h(v)=\left\{
\begin{array}{ll}
1,&\textup{ if }v=u,\\
g(v),&\textup{ if }v\neq u
\end{array}\right.
$$
is a resolving function of $G$. Hence $\dim_f(G)\leq\rdim_f(G^u)+1$,
as desired. $\qed$

If $u$ is not an end-vertex of the path $P_n$, then
$\rdim_f(P_n^u)=\rdim(P_n^u)=\dim_f(P_n)=1$, which implies that the
upper bounds in Proposition~\ref{root frac} (i) and (iii) are tight.
The fact that $\rdim_f(K_n^u)=\frac{n-1}{2}$ shows that the
inequality in Proposition~\ref{root frac} (ii) is tight.

\medskip
Next, we study the fractional metric dimension for the hierarchical product of
graphs.

For two rooted graphs $G_1^{u_1}$ and $G_2^{u_2}$, write
$$
\begin{array}{c}
\mathcal P^{u_1}=\{\{x,y\}\subseteq V(G_1)\mid
x\neq
y,d_{G_1}(x,u_1)=d_{G_1}(y,u_1)\},\\
\overline{\mathcal P}^{u_2u_1}=\{\{x_2x_1,y_2y_1\}\subseteq V(G_2^{u_2}\sqcap G_1^{u_1})\mid x_2x_1\neq y_2y_1, d(x_2x_1,u_2u_1)=d(y_2y_1,u_2u_1)\}.
\end{array}
$$

\begin{lemma}\label{fraclow}
Let $G_1^{u_1}$ and $G_2^{u_2}$ be two rooted graphs.  If $G_1$ is
not a path with an end-vertex $u_1$, then
$$
\rdim_f(G_2^{u_2}\sqcap G_1^{u_1})\geq
|V(G_2)|\cdot\rdim_f(G_1^{u_1}).
$$
\end{lemma}
\proof Suppose that $\overline g$ is a rooted resolving function of
$G_2^{u_2}\sqcap G_1^{u_1}$ with weight $\rdim_f(G_2^{u_2}\sqcap G_1^{u_1})$. For each $z\in V(G_2)$,
define
$$
\overline g_z: V(G_1)\longrightarrow[0,1],\quad x\longmapsto\overline g(zx).
$$
Write $\overline{\mathcal P}^{u_1}=\{\{zx,zy\}\mid z\in
V(G_2),\{x,y\}\in\mathcal P^{u_1}\}$. By (\ref{1}), we have
$\overline{\mathcal P}^{u_1}\subseteq\overline{\mathcal
P}^{u_2u_1}.$ Hence $\overline g_z(R_{G_1}\{x,y\})\geq 1$ for any
$\{x,y\}\in\mathcal P^{u_1}$, which implies that $|\overline
g_z|\geq\rdim_f(G_1^{u_1})$. Consequently,
$$
\rdim_f(G_2^{u_2}\sqcap G_1^{u_1})=|\overline g|=\sum_{z\in
V(G_2)}|\overline g_z|\geq |V(G_2)|\cdot\rdim_f(G_1^{u_1}),
$$
as desired.
$\qed$

\begin{thm}\label{fracmetmain}
Let $G_1^{u_1}$ and $G_2^{u_2}$ be two rooted graphs. If $G_1$ is
not a path with an end-vertex $u_1$, then
$$
\dim(G_2^{u_2}\sqcap G_1^{u_1})=|V(G_2)|\cdot\rdim_f(G_1^{u_1}).
$$
\end{thm}
\proof Combining Proposition \ref{root frac} and Lemma
\ref{fraclow}, we only need to prove that
\begin{equation}\label{10}
\dim_f(G_2^{u_2}\sqcap
G_1^{u_1})\leq|V(G_2)|\cdot\rdim_f(G_1^{u_1}).
\end{equation}
By Proposition \ref{path} we have $\mathcal P^{u_1}\neq\emptyset$.
Let $g$ be a rooted resolving function of $G_1$ with weight
$\rdim_f(G_1^{u_1})$. Define
$$
\overline g: V(G_2^{u_2}\sqcap G_1^{u_1})\longrightarrow[0,1],\quad x_2x_1\longmapsto g(x_1).
$$
We shall show that, for any two distinct vertices $x_2x_1$ and $y_2y_1$ of $G_2^{u_2}\sqcap G_1^{u_1}$,
\begin{equation}\label{11}
\overline g(R\{x_2x_1,y_2y_1\})\geq 1.
\end{equation}

{\em Case 1.} $x_2=y_2.$ If $u_1\not\in R_{G_1}\{x_1,y_1\}$, by
Lemma \ref{resolve} we get $R\{x_2x_1,y_2y_1\}=\{x_2z\mid z\in
R_{G_1}\{x_1,y_1\}\}$, which implies that $\overline
g(R\{x_2x_1,y_2y_1\})= g(R_{G_1}\{x_1,y_1\})$. Since
$\{x_1,y_1\}\in\mathcal P^{u_1}$, we obtain (\ref{11}). If $u_1\in
R_{G_1}\{x_1,y_1\}$, by Lemma \ref{resolve} we have
$R\{x_2x_1,y_2y_1\}\supseteq\{vz\mid z\in V(G_1)\}$ for any $v\in
V(G_2)\setminus\{x_2\}$, which implies that $\overline
g(R\{x_2x_1,y_2y_1\})\geq |g|$, so (\ref{11}) holds.

{\em Case 2.} $x_2\neq y_2.$ Write $W=\{z\mid x_2z\in
R\{x_2x_1,y_2y_1\}\}$ and $S=\{z\mid y_2z\in R\{x_2x_1,y_2y_1\}\}$.
By Lemma \ref{resolve} we have $W\cup S=V(G_1)$. Then
$$
\overline g(R\{x_2x_1,y_2y_1\})\geq\sum_{z\in W}\overline g(x_2z)+\sum_{z\in S}
=g(W)+g(S)\overline g(y_2z)\geq |g|,
$$
which implies that (\ref{11}) holds.

Therefore, $\overline g$ is a resolving function of $G_2^{u_2}\sqcap
G_1^{u_1}$, which implies that $\dim_f(G_2^{u_2}\sqcap
G_1^{u_1})\leq|\overline g|$. Since $|\overline
g|=|V(G_2)|\cdot\rdim_f(G_1^{u_1}),$ we obtain (\ref{10}). Our proof
is accomplished. $\qed$

\begin{cor}\label{34}
Let $n\geq 2$. Then $\dim_f(T_n)=2^{n-2}$.
\end{cor}
\proof It is immediate from Theorem~\ref{fracmetmain}. $\qed$

By Corollaries~\ref{29} and \ref{34}, the binomial tree $T_n$ is a
graph whose metric dimension is equal to its fractional metric
dimension.

Finally, we shall prove some tight inequalities for
$\dim_f(G^{u}\sqcap P_n^0)$.

\begin{prop}\label{frapathbound}
For any rooted graph $G^u$, we have
$$
\dim_f(G)\leq\dim_f(G^u\sqcap P_n^0)\leq\dim_f(G^u\sqcap P_{n+1}^0)\leq\frac{|V(G)|}{2}.
$$
\end{prop}
\proof Write $P_n=(z_0=0,z_1,\ldots,z_{n-1})$. For a resolving
function $\overline g_{n+1}$ of $G^u\sqcap P_{n+1}^0$, we define
$\overline g_{n+1}': V(G^u\sqcap P_n^0)\longrightarrow[0,1]$ by
$$
\overline g_{n+1}'(x_2x_1)=\left\{
\begin{array}{ll}
\overline g_{n+1}(x_2z_{n-1})+\overline g_{n+1}(x_2z_n),& \textup{ if }x_1=z_{n-1},\\
\overline g_{n+1}(x_2x_1), &\textup{ if }x_1\neq z_{n-1}.
\end{array}\right.
$$
Then $\overline g_{n+1}'$ is a resolving function of $G^u\sqcap
P_n^0$. Since $|\overline g_{n+1}'|=|\overline g_{n+1}|$, we have
$$
\dim_f(G)=\dim_f(G^u\sqcap P_1^0)\leq\dim_f(G^u\sqcap
P_n^0)\leq\dim_f(G^u\sqcap P_{n+1}^0).
$$

Define $\overline h: V(G^u\sqcap P_{n+1}^0)\longrightarrow[0,1]$ by
$$
\overline h(x_2x_1)=\left\{
\begin{array}{ll}
\frac{1}{2},& \textup{ if }x_1=z_{n},\\
0, &\textup{ if }x_1\neq z_n.
\end{array}\right.
$$
Then $\overline h$ is a resolving function of $G^u\sqcap P_{n+1}^0$
with weight $\frac{|V(G)|}{2}.$ Hence
$\dim_f(G^u\sqcap P_{n+1}^0)\leq\frac{|V(G)|}{2}$.
$\qed$

For $m\geq 2$, we have $\dim_f(K_m^u\sqcap P_n^0)=\frac{m}{2}$. This
shows that all the inequalities in Proposition~\ref{frapathbound}
are tight.

\section*{Acknowledgement}
 This research is supported by NSFC(11271047),
SRFDP and the Fundamental Research Funds for the Central University
of China.

\end{document}